\documentclass[11pt,leqno]{amsart}
\usepackage{minart,amsmath,amssymb}
\title{Quaternion Orders and Ternary Quadratic Forms}
\author{Stefan Lemurell}

\newcommand{\rr}{\ensuremath{\mathbb R}}

\newcommand{\q}{\ensuremath{\mathbb Q}}
\newcommand{\h}{\ensuremath{\mathbb H}}

\newcommand{\fa}{\ensuremath{\mathfrak a}}

\newcommand{\id}{\ensuremath{\mathfrak I}}
\newcommand{\p}{\ensuremath{\mathfrak p}}
\newcommand{\A}{\ensuremath{\mathfrak A}}

\newcommand{\OO}{\ensuremath{\mathcal O}}

\newcommand{\md}[2]{\equiv #1\;(\mbox{mod } #2)}

\newcommand{\Matrix}[4]{\left( \begin{array}{cc} #1 & #2 \\#3 &
    #4 \end{array} \right)}  
\newcommand{\ter}[3]{\left<#1\right>\perp\left<#2\right>\perp\left<#3\right>}  

\begin{document}
\bibliographystyle{siam}

\maketitle

\section*{Introduction}\label{sec:int}
The main purpose of this paper is to provide an introduction to
the arithmetic theory of quaternion algebras. However, it also
contains some new 
results, most notably in Section~\ref{sec:loc}. We will emphasise on
the connection between quaternion algebras and quadratic forms. This
connection will provide us with an efficient tool to consider
arbitrary orders instead of having to restrict to special classes of
them. The existing results are mostly restricted to special classes of
orders, most notably to so called Eichler orders.

The paper is organised as follows.  Some
notations and background are provided in Section~\ref{sec:not},
especially on the theory of quadratic forms. Section~\ref{sec:alg}
contains the basic theory of quaternion algebras. Moreover at the end
of that section, we give a quite general solution to the problem of
representing a quaternion algebra with given discriminant. Such a
general description seems to be lacking in the literature.

Section~\ref{sec:ord} gives the basic definitions concerning orders in
quaternion algebras. In Section~\ref{sec:ter}, we prove an important
correspondence between ternary quadratic forms and quaternion
orders. Section~\ref{sec:loc} deals with orders in quaternion algebras
over $\p$-adic fields. The major part is an investigation of the
isomorphism classes in the non-dyadic and $2$-adic cases. The
starting-point is the correspondence with ternary quadratic forms and
known classifications of such forms. From this, we derive
representatives of the isomorphism classes of quaternion orders. These
new results are complements to existing more ring-theoretic descriptions of
orders. In particular, they are useful for computations.

Finally, section~\ref{sec:glo} contains the basic theory of orders in
quaternion algebras over algebraic number fields and the connection
with the $\p$-adic case. At the end of that section, we give an
explicit basis of a maximal order when the discriminant of the algebra
is a principal ideal. It is related to the results at the end of
Section~\ref{sec:alg} and is more general than any other similar
description known to the author.

\section{Notation and some background} \label{sec:not}
For convenience of the reader, we will recall
the basic definitions of the theory of quadratic forms. Let $K$ be a
field with 
char$K\neq2$ and $R$ any subring of $K$ containing $1$.  To every
quadratic form 
\begin{equation}
  \label{form}
  f=f(X_{1},\ldots,X_{n})=\sum_{1\leq i\leq j\leq n}a_{ij}X_{i}X_{j},   
\end{equation}
with coefficients $a_{ij}\in K$, we can associate a symmetric $n\times
n$-matrix $M_{f}=(m_{ij})$, where 
\[  m_{ij}=\left\{ 
  \begin{array}{ll}
    a_{ij}, & i<j \\
    2a_{ii}, & i=j \\
    a_{ji}, & i>j,
  \end{array}
\right. \]
and the corresponding bilinear form $B_{f}(x,y)=x^{t}M_{f}y$. A direct
inspection gives $f(x)=\frac12 B_{f}(x,x)$. This is where we get into
trouble, if char$K=2$. We will only consider non-degenerate forms, that is,
forms for which \mbox{$\det(M_{f})\neq0$}. The form $f$ is said to represent
$a\in K$ over $R$, if there exists $x\in R^{n}$ such that $f(x)=a$. It
is called isotropic, if there is a nontrivial representation of $0$,
otherwise it is called anisotropic. We will use the standard notation
$f=\left< a_{1}\right> \perp\cdots\perp\left< a_{n}\right>$ to denote
the diagonal form $f= a_{1}X_{1}^{2}+\cdots+ a_{n}X_{n}^{2}$.

Two forms $f$ and $g$ are called
isometric over $R$, $f\cong g$, if  there exists $T\in GL_{n}(R)$ such
that $M_{f}=T^{t}M_{g}T$. They are said to be similar over $R$,
denoted by $f\sim g$, if
there exists an $u\in R^{*}$ such that $u\cdot f$ is isometric to $g$. Both
isometry and similarity are equivalence relations, similarity
obviously being coarser.  

A form $f$ like \mref{form} is called integral over $R$ if $a_{ij}\in
R$. It is called primitive, if  the ideal generated by the
coefficients $a_{ij}$ is equal to $R$. 

The discriminant, $d(f)$, of a non-degenerate quadratic form $f$ is
defined to be 
the class of $\det(M_{f})$ in $K^{*}/(K^{*})^{2}$. The reason for
taking classes modulo $(K^{*})^{2}$ is to make it an invariant of
isometry classes. Note that it is only an invariant of similarity
classes in even dimensions, since multiplication by $u$ multiplies the
determinant by $u^{n}$. If the form $f$ is integral over $R$, then
the discriminant of $f$ regarded as a form over $R$ is the class of
$\det(M_{f})$ as an element in the multiplicative set
$R\setminus\{0\}$ modulo $(R^{*})^{2}$. In the case of 
quadratic forms of odd dimensions, it is customary and natural to take
the discriminant 
to be the class of $\frac12\det(M_{f})$ instead, and we will follow this
convention.

Let $K$ be an algebraic number field with ring of integers $R$. This
will be the case for much of the paper and exceptions will be clearly
notified. Let $\Omega$ be the set of places (normalised valuations) on
$K$, $\Omega_{f}$ the finite (non-archimedean) and
$\Omega_{\infty}$  the infinite (archimedean) ones. If $\nu\in\Omega$,
then $K_{\nu}$ will denote the completion of $K$ with respect to
$\nu$, and if $\nu\in\Omega_{f}$, then $R_{\nu}$ will be the ring of
integers in $K_{\nu}$. 

When a ring $A$ is understood, then $\left< x_{1},\ldots,x_{n}\right>$
will denote the free $A$-module generated by $\{x_{1},\ldots,x_{n}\}$.

\section{Quaternion algebras} \label{sec:alg}
In this section, we will give the definition and some fundamental
properties of quaternion algebras. We will reduce in generality along
the way to make the exposition as simple and clear as possible.

\defn{Let $K$ be an arbitrary field. A quaternion algebra $\A$ over
  $K$ is a simple, central algebra of dimension $4$ over $K$.}

From Wedderburn's structure theorem on simple algebras
\cite[2.5]{Kersten}, one deduces that either $\A\cong M_{2}(K)$ or
$\A\cong D$, where $D$ is a division algebra with centre $K$. 

We will
from now on assume that char$K\neq 2$. With this condition, it is
always possible to find a convenient `diagonal' basis $1,i,j,ij$ of
$\A$ over $K$, which satisfies
\begin{equation} \label{algbasis} 
i^{2}=a,\;j^{2}=b \mbox{ and } ij=-ji, \mbox{ where } a,b\in
K,\;ab\neq 0. 
\end{equation}
A proof of this can be found in a more general setting in \cite[\S
17]{Kersten}. We will denote such an algebra by $(a,b)_{K}$. 

\remn{
To make everything completely explicit, we remark that it is possible
to embed $(a,b)_{K}$ in $M_{2}(K(\sqrt{a}))$  for
example by
\[ i\longmapsto\Matrix{\sqrt{a}}{0}{0}{-\sqrt{a}},\;
j\longmapsto\Matrix{0}{b}{1}{0}. \]
From this it is clear that if $a$ is a square in $K$, then
$(a,b)_{K}\cong M_{2}(K)$. A necessary and sufficient condition for
$(a,b)_{K}\cong M_{2}(K)$ is that $a$ is the norm of an element in
$K(\sqrt{b})$ with respect to $K$
\cite[17.4]{Kersten}. Of course, one may interchange $a$ and $b$ in
this remark.
}

There is a natural involution in $\A$, which in a basis satisfying
\mref{algbasis} is given by
\[ x=x_{0}+x_{1}i+x_{2}j+x_{3}ij \longmapsto \bar{x}=
x_{0}-x_{1}i-x_{2}j-x_{3}ij. \]
One defines the (reduced) trace and (reduced)
norm from $\A$ into $K$ by
\[ Tr(x)=x+\bar{x} \mbox{ and } N(x)=x\bar{x}. \]
The norm is a quaternary quadratic form over $K$, with 
corresponding symmetric bilinear form given by
$B(x,y)=Tr(x\bar{y})$. A direct calculation gives 
\[ N(x)=x_{0}^{2}-ax_{1}^{2}-bx_{2}^{2}+abx_{3}^{2}, \]
if $x=x_{0}+x_{1}i+x_{2}j+x_{3}ij\in (a,b)_{K}$. From this, we see that
the determinant of the norm on $(a,b)_{K}$ is equal to
$16a^{2}b^{2}$. Hence, the discriminant of the norm form of a
quaternion algebra is equal to $1$.

Now define the set of pure
quaternions $\A_{0}$ to be
\[ \A_{0}=\{q\in\A : Tr(q)=0\}, \]
and denote the norm $N$ restricted to $\A_{0}$ by $N_{0}$. If we have
chosen a basis of $\A$ satisfying \mref{algbasis}, then obviously
$\A_{0}=\left<i,j,ij\right>$. We conclude that $\A_{0}$ is a
\mbox{$3$-dimensional} $K$-vector space and $N_{0}$ a ternary quadratic
form. A nontrivial and interesting fact is that $N$ is isotropic iff
$N_{0}$ is isotropic \cite[42:12]{O'Meara}.

The norm form determines the quaternion algebra in the following sense
\cite[\S 57]{O'Meara}:
\prop{ \label{equivalent}
Let $K$ be a field with char$K\neq 2$ and let $\A$ and 
$\A^{'}$ be  quaternion algebras over $K$ with
corresponding norm forms $N$ and $N^{'}$. Then the following statements
are equivalent:
\begin{enumerate}
\item $\A \cong\A^{'}$
\item $N$ and $N^{'}$ are isometric
\item $N_{0}$ and $N_{0}^{'}$ are isometric
\end{enumerate}
Moreover, $\A$ a division algebra is equivalent to $N$ being
anisotropic, which in turn is equivalent to  $N_{0}$ being anisotropic.
}

This correspondence between quaternion algebras and quaternary and
ter\-nary quadratic forms will be used in the sequel. Especially a
refinement of the correspondence with ternary quadratic forms will be
presented in detail in Section~\ref{sec:ter} and used in the following
sections.

We now specialise to the fields of our main interest. So from now on
assume that $K$ is an algebraic number field. If $\nu$ is a place on
$K$, then $\A_{\nu}=\A\otimes_{K}K_{\nu}$ is a quaternion algebra over
$K_{\nu}$. 

Every quaternary quadratic form $Q$ over $K$ represents $1$ iff it is not
negative definite at any real place \cite[Satz 19]{Hasse6}. Hence, we
may conclude that 
\[ Q\cong \left< 1 \right>\perp\left< a \right>\perp\left< b
\right>\perp\left< c  \right>, \]
for some $a,b,c\in K$, if $Q$ is not
negative definite at any real place. If the discriminant of $Q$ is
equal to $1$, then we may assume that $c=ab$ and we conclude from
\mref{equivalent} the following: 

\prop{ \label{corr}
Let $K$ be an algebraic number field. Then there is a one-to-one
correspondence between 

\begin{enumerate}
\item isomorphism classes of quaternion algebras
over $K$,
\item isometry classes of quaternary quadratic forms over
$K$ with discriminant  equal to $1$, which are not
negative definite at any real place.
\end{enumerate}

\noindent Furthermore, if $\nu$ is a
place on $K$, then there is a one-to-one
correspondence between

\begin{enumerate}
\item isomorphism classes of quaternion algebras
over $K_{\nu}$, 
\item isometry classes of quaternary quadratic forms over
$K_{\nu}$ with discriminant equal to $1$, which are not
negative definite if $K_{\nu}=\rr$.
\end{enumerate}
}

Hence, one can use the classification of quadratic forms in order to
classify quaternion algebras up to isomorphism.

The classification of quadratic forms up to isometry over $\q$ was
first made by Minkowski in \cite{Minkowski}. This classification was
later simplified and generalised to arbitrary algebraic number fields
by Hasse in \cite{Hasse2} and \cite{Hasse7}. Using the results of
Hasse and \mref{corr} one derives the following two theorems on the
classification of quaternion algebras.

\thm{
Let $K_{\nu}$ be a completion of a number field $K$. Then there are exactly
two quaternion algebras over $K_{\nu}$ up to isomorphism, $M_{2}(K_{\nu})$
and a division algebra. 
}

The division algebra over $K_{\nu}$ will be denoted by $\h_{\nu}$. We
 say that $\A$ is ramified  at $\nu$ if $\A_{\nu}\cong \h_{\nu}$, 
otherwise $\A$ is said to split at $\nu$. If $\nu$ is a real place,
then one often uses definite/indefinite instead of
ramified/split, especially if $K=\q$.

\thm{ \label{classification}
Let $K$ be an algebraic number field and $\A$ and $\A^{'}$ two
quaternion algebras over $K$. Then the following statements are
equivalent:
\begin{enumerate}
  \item $\A\cong\A^{'}$,
  \item $\A_{\nu}\cong\A^{'}_{\nu},\;\forall \nu\in\Omega$,
  \item $\A$ and $\A^{'}$ are ramified at the same places.
\end{enumerate}
Moreover, $\A$ is always ramified at an even number of
places. Conversely, given an even number of places, it is always
possible to find a quaternion algebra which is ramified at exactly
these places.
}

The (reduced) discriminant, $d(\A)$, of a quaternion algebra \A\ is
defined to 
be the product of the prime ideals $\p$ at which $\A$ is ramified. This
is a well-defined invariant of the isomorphism classes by
\mref{classification}. If $K=\q$, then the 
discriminant determines the isomorphism class, but for other fields
one clearly also needs information on the infinite ramifications. 

With this classification at hand, it is of course of interest to be
able to easily determine the ramifications of a given quaternion
algebra over $K$. It turns out that the  Hilbert symbol solves the
problem. If $\nu\in\Omega$ and $a,b\in K_{\nu}^{*}$,
then the Hilbert symbol $(a,b)_{\nu}$ in $K_{\nu}$ is defined by 
\[ (a,b)_{\nu}= \left\{ \begin{array}{rl}
 1,  & \mbox{if }  X^{2}-aY^{2}-bZ^{2} \mbox{ is isotropic over } K_{\nu} \\
-1,  &   \mbox{otherwise.}         
\end{array} \right. \]
We remark that if $N_{0}$ is the restricted norm form on $(a,b)_{K}$, then
\[ abN_{0}\cong\ter{-a}{-b}{1}. \]
From this we conclude that $N_{0}$ is isotropic at $\nu$ iff
$(a,b)_{\nu}=1$. Hence, it is immediate from \mref{equivalent} that
$(a,b)_{K}$ is ramified at $\nu$ iff $(a,b)_{\nu}=-1$.
For a proof of the following properties of the Hilbert symbol, see
\cite[\S 63,71]{O'Meara}. We remark, that property (d) is exactly what
is needed to prove that $\A$ is always ramified at an even number of
places. 

\prop{ \label{hilb}
The Hilbert symbol satisfies the following:
\newline 
  $ \begin{array}{ll}
       &  \\
    \mbox{\rm (a) } & (a,bc)_{\nu}=(a,b)_{\nu}(a,c)_{\nu},\; (a,-a)_{\nu}=1
                \mbox{ and } (a,b^{2})_{\nu}=1. \\
    \mbox{\rm (b) } & \mbox{If $\nu$ is real, then } (a,b)_{\nu}=-1
    \mbox{ iff } a<0 \mbox{ and } b<0. \\ 
    \mbox{\rm (c) } & \mbox{If $\p$ is a non-dyadic prime, then }
        \left\{ \begin{array}{ll}
        (a,b)_{\p}=1,  &  \mbox{if } a,b \in R_{\p}^{*} \\
        (a,\p)_{\p}=(\frac{a}{\p}), & \mbox{if } a \in R_{\p}^{*},
      \end{array} \right. \\
                & \mbox{where } (\frac{a}{\p})=1  \mbox{ if $a$ is a
                   square modulo } \p \mbox{ and $-1$ otherwise.} \\
     \mbox{\rm (d) } & \prod_{\nu \in \Omega} (a,b)_{\nu}=1,
      \; \; \forall a,b \in K^{*}. 
  \end{array} $
}

The dyadic case is more difficult in general. If there is only
one dyadic prime, then \mref{hilb}(d) solves the problem, but otherwise
we might have to do some calculations. However, it is possible to
prove the following
general result. If $\pi$ is a prime element in
$R_{\p}$, $\epsilon\in R_{\p}^{*}$ and $\alpha\in R_{\p}^{*}$ is an
element of 
quadratic defect $\delta_{\p}(\alpha)=4R_{\p}$, then \cite[63:11a]{O'Meara}
\begin{equation} \label{defect}
(\pi,\alpha)_{\p}=-1 \mbox{ and }  (\epsilon,\alpha)_{\p}=1.
\end{equation}
(An element is of quadratic defect $4R_{\p}$, if it is not a square but
congruent to a square modulo $4R_{\p}$.)

The concluding result of this section gives an answer to the question:
Given the ramifications of a quaternion algebra $\A$ over $K$, how to find
\mbox{$a,b\in R$} such that $\A\cong (a,b)_{K}$? If
$d(\A)=\p_{1}\cdots \p_{r}$, then let $n_{i}$ be  positive
integers such that $\id=\p_{1}^{n_{1}}\cdots \p_{r}^{n_{r}}$ is
principal. A generator of the principal ideal $\id$  will be
called a representative of $d(\A)$. The proposition below will
give an answer to our question when $d(\A)=R$ or $d(\A)$ has a representative
$\p_{1}^{n_{1}}\cdots \p_{r}^{n_{r}}$, such that all integers $n_{i}$
are odd,  in particular, when $d(\A)$ is principal. It is also a very 
explicit proof of the 
last conclusion in \mref{classification} in this restricted case. In
principal, it follows the same 
idea as in Hasse's original proof.

\prop{ \label{explicit}
Let $\A$ be a quaternion algebra over an algebraic number field
$K$. Suppose that $d(\A)=R$ or $d(\A)=\p_{1}\cdots \p_{r}$ has a representative
$(d)=\p_{1}^{n_{1}}\cdots \p_{r}^{n_{r}}$, such that all integers $n_{i}$
are odd. Then choose $a\in R$ to be a generator of a prime ideal satisfying 
{\em gcd}$(a,d)=1$ and
\begin{equation} \label{alpha}
\left\{ 
  \begin{tabular}{ll}
    $a<0$, & at $\nu$ real, if $\A$ ramified at $\nu$, \\
    $a>0$, & at $\nu$ real, if $\A$ split at $\nu$, \\
    $(\frac{a}{\p})=-1$, & at $\p$ non-dyadic, if $\A$ ramified at
    $\p$, \\
    $\delta_{\p}(a)=4R_{\p}$, & at all dyadic primes $\p$. 
  \end{tabular} \right.
\end{equation}
Then $\A\cong (a,-d)_{K}$.
}
\begin{proof}
  The existence of such an $a$ is assured by the generalisation to arbitrary
  number fields of
  Dirichlet's Theorem on primes in linear progressions 
  \cite[Cor.4, p.360]{Narkiewicz}. 
  The ramifications of $\A$ and $(a,-d)_{K}$ obviously agree at
  all real places. 
  If $\p$ is a non-dyadic prime dividing $d$, then
  \mref{hilb} implies that 
  \[(a,-d)_{\p}=(a,\p)_{\p}=(\frac{a}{\p})=-1.\]
  This is where it is necessary to have $\p$ to an odd power in $(d)$.
  The only other non-dyadic prime at which $(a,-d)_{K}$ could
  possibly be ramified is $a$, since $a,d\in R_{\p}^{*}$ for
  all other non-dyadic primes. 

  If $\p$ is a dyadic prime, then \mref{defect}
  implies that $(a,-d)_{\p}=-1$ iff $\p|d$. Hence, the
  ramifications of $\A$ and $(a,-d)_{K}$ differ at most at
  $a$ and then \mref{classification} implies that $\A\cong
  (a,-d)_{K}$. 
\end{proof}

For the field of rational numbers \mref{explicit} simplifies to
(compare \cite{Latimer})
\cor{
  Let $\A$ be an indefinite quaternion algebra over $\q$
with discriminant $d =p_{1} \cdots p_{2r}$. Choose $p$ to
be a prime such that $p \md58$ and 
$(\frac{p}{p_{i}}) =-1, \; \forall p_{i}>2$.
 Then ${\A} \cong (p,d)_{\q} \cong (p,-d)_{\q}$.
}
\cor{
  Let $\A$ be a definite quaternion algebra over $\q$
with discriminant $d =p_{1} \cdots p_{2r-1}$. Choose $p$ to
be a prime such that $p \md38$ and 
$(\frac{p}{p_{i}}) =-1, \; \forall p_{i}>2$.
 Then ${\A} \cong (-p,-d)_{\q}$.
}

\section{Quaternion orders} \label{sec:ord}
In this section, $R$  will be a Dedekind ring with field of fractions $K$. 

An $R$-lattice $\Lambda$ on $\A$ is a finitely generated $R$-module
such that $K\Lambda=\A$. We remark that it is not always
possible to find an 
$R$-basis of $\Lambda$. However, it is always possible to find a basis
$e_{1},\ldots,e_{4}$ of 
$\A$ and an $R$-ideal $\fa$, such that
$\Lambda=\fa e_{1}\oplus R e_{2}\oplus R e_{3}\oplus R e_{4}$
\cite[81:5]{O'Meara}.  Let $\Lambda$ and $\Lambda^{'}$
be two lattices on \A.
The index,  $[\Lambda :\Lambda^{'}]$, of $\Lambda^{'}$ in $\Lambda$ is
defined to be the $R$-ideal generated by $\det(\varphi)$ of all linear
transformations 
$\varphi : \A\longrightarrow\A, \mbox{ such that }
\varphi(\Lambda)\subseteq \Lambda^{'}. $
In particular,
if $\Lambda$ and $\Lambda^{'}$ both have $R$-bases, then $[\Lambda
:\Lambda^{'}]$ is the determinant of the matrix which takes a basis of
$\Lambda$ into a basis of $\Lambda^{'}$. Given an $R$-lattice
$\Lambda$ on \A, we define its dual $\Lambda^{\#}$ to be
\[ \Lambda^{\#}=\{ q\in\A : Tr(q\Lambda)\subseteq R\}. \]
This is again an $R$-lattice on $\A$.

\defn{
An ($R$-)order $\OO$ in a quaternion algebra $\A$ is an
\mbox{$R$-lattice on $\A$,} 
which is also a ring containing $R$.
}

For the rest of this paper $\OO$ will always denote a quaternion
order. We remark that if $x\in\OO$, then $x$ is integral over $R$, so
$Tr(x),N(x)\in R$.

If $\Lambda$ is a lattice on $\A$, then we define the left (right)
order $\OO_{l}(\Lambda)$ ($\OO_{r}(\Lambda)$) of $\Lambda$ to be
\[ \OO_{l}(\Lambda)=\{x\in\A : x\Lambda\subseteq\Lambda\} \mbox{ and }
\OO_{r}(\Lambda)=\{x\in\A : \Lambda x\subseteq\Lambda\}. \]
It is easy to verify that both $\OO_{l}(\Lambda)$ and $\OO_{r}(\Lambda)$
are orders in \A.

A left (right) ideal of an order \OO\ is a lattice $\Lambda$ on $\A$ such
that $\OO\Lambda\subseteq\Lambda$ ($\Lambda\OO\subseteq\Lambda$). The
ideal is called a two-sided ideal of \OO, if it is both a left and
right ideal of $\OO$. If $\Lambda$ is any lattice, then obviously $\Lambda$ is a
left $\OO_{l}(\Lambda)$-ideal and a right $\OO_{r}(\Lambda)$-ideal. An
ideal $\Lambda$ is a principal \OO-ideal, if there is $x\in\A$ such that
$\Lambda=\OO x$.

The most important invariant of a quaternion order is the (reduced)
discriminant. It is defined as follows. Let $\id$ be the ideal
generated by all $\det(Tr(x_{i}\bar{x}_{j}))$, where $x_{1},\ldots,
x_{4}\in\OO$. It is easy to prove that $\id$ is the square of an ideal.
The (reduced) discriminant, $d(\OO)$, of $\OO$ is defined to be the
square root of $\id$. If $\OO$ and $\OO^{'}$ are two orders, then the
discriminants satisfy
\begin{equation}
  \label{discindex}
  d(\OO)=d(\OO^{'})\cdot [\OO^{'} : \OO]. 
\end{equation}
\exn{
Let $\A=(a,b)_{K}$ be a quaternion algebra. Suppose that $a,b\in
R$. This is no restriction, since $(a,b)_{K}\cong (ax^{2},by^{2})_{K}$
for all $x,y\in K$. Then 
\[ \OO=R+Ri+Rj+Rij, \]
where $i,j$ satisfy \mref{algbasis} is an order in \A. We will
denote this order by $(a,b)_{R}$. The discriminant of $(a,b)_{R}$ is
\[ (\det(Tr(x_{i}\bar{x}_{j})))^{\frac12},\] 
where
$\{x_{1},\ldots,x_{4}\}=\{1,\,i,\,j,\,ij\}$. A simple calculation
shows that the matrix is diagonal
and  $d(\OO)=(4ab)$.
}

A maximal order is an order, which is not strictly contained in any other
order. The concept of maximal orders is more complicated than in the
commutative case, since there is more than one maximal order in general.
This complication occurs since $\{x \in\A: \: N(x),Tr(x) \in R \}$ 
is not always a ring. However, the discriminants of all maximal orders
in a quaternion algebra always agree (see \mref{maximal} and
\mref{maxcond}). 

The task to classify quaternion orders is of course much more
complicated than for quaternion algebras. In particular, there is no
analogue of \mref{classification}, that is, isomorphism at all places
no longer implies global isomorphism. However, the investigation of
$\OO_{\nu}=\OO\otimes_{R}R_{\nu}$ is still essential and useful in order
to classify quaternion orders $\OO$ in algebras over algebraic number
fields (see Section~\ref{sec:glo}). 

In order to get a structure on the set of orders in a
quaternion algebra, we introduce special classes of orders. An order 
$\OO$ is called a Gorenstein order if $\OO^{\#}$
is projective as $\OO$-module, and it is called a Bass order if every order
in $\A$ containing $\OO$ is a Gorenstein order. For an arbitrary
order $\OO$, there is a unique Gorenstein order $G(\OO)$ (called the 
Gorenstein closure) and a unique $R$-ideal $b(\OO) \subseteq R$ (the
Brandt invariant) such that \cite[(1.4)]{Brzezinski2}
\begin{equation}
  \label{gorenstein}
  \OO=R+b(\OO)G(\OO).
\end{equation}

In \cite{Eichler2}, Eichler introduced so called primitive orders. An
order is called primitive, if it contains a
maximal order of a quadratic subfield of $\A$. It is easy to show
that every primitive order is a Bass order. However, whether every
Bass order is primitive is an open question in general. It is true in
the local case by \cite[(1.11)]{Brzezinski4} and in the case of rational
orders by \cite[Satz 8]{Eichler2}.

An order $\OO$ is called a hereditary order if every $\OO$-ideal is
$\OO$-projective. It is well-known that hereditary orders are exactly
those with square free discriminant \cite[(39.14)]{Reiner}. These
classes of orders obviously satisfy the inclusions
\[ \mbox{\{Gorenstein\}}\supset \mbox{\{Bass\}}\supset
\mbox{\{Hereditary\}}\supseteq \mbox{\{Maximal\}}. \]

\section{Ternary quadratic forms and quaternion orders}
\label{sec:ter} 
Before we continue the investigation of quaternion orders, we will in
this section refine the correspondence between quaternions and ternary
quadratic forms. We will have to assume that $R$ is a principal ideal
domain, since otherwise we will have some complications.

If $f$ is a non-degenerate ternary quadratic form integral over $R$, then
define $C_{0}(f)$ to be the even Clifford algebra over $R$ 
associated to $f$. Then $C_{0}(f)$ is an order in a quaternion algebra
over $K$. If 
\[ f=\sum_{1 \leq i \leq j \leq 3} a_{ij}X_{i}X_{j}, \]
then a direct computation shows that $C_{0}(f)$ has an  $R$-basis
$1,e_{1},e_{2},e_{3}$ such that: 
\begin{eqnarray} \label{cliffrel}
  e_{i}^{2} & = & a_{jk}e_{i}-a_{jj}a_{kk}, \nonumber \\
  e_{i}e_{j} & = & a_{kk}(a_{ij}-e_{k}), \\
  e_{j}e_{i} & = & a_{1k}e_{1}+a_{2k}e_{2}+a_{3k}e_{3}-a_{ik}a_{jk}, \nonumber
\end{eqnarray}
where $(i,j,k)$ is an even permutation of $(1,2,3)$. From this we get
that the norm form for $C_{0}(f)$ in this basis is
\begin{equation}
  \label{norm}
  Q=X_{0}^{2}+\sum_{(i,j,k)}\left[a_{ij}X_{0}X_{k}+
    a_{ii}a_{jj}X_{k}^{2}+(a_{ik}a_{jk}-a_{ij}a_{kk})X_{i}X_{j}\right], 
\end{equation}
where the sum is over all even permutations $(i,j,k)$ of $(1,2,3)$. 
It is trivial to
check from the relations \mref{cliffrel} that if $\epsilon\in R^{*}$, then
$C_{0}(f)=C_{0}(\epsilon f)$.

Conversely, assume that $\OO=\left<1,e_{1},e_{2},e_{3} \right>$ is an
$R$-order in a quaternion algebra 
$\A$. Then
$\Lambda=\OO^{\#}\cap \A_{0}$ is a $3$-dimensional $R$-lattice on
$\A_{0}$. If $\OO^{\#}=\left<f_{0},f_{1},f_{2},f_{3}\right>$ where
$\{f_{i}\}$ is the dual basis of $\{e_{i}\}$, then
$\Lambda=\left<f_{1},f_{2},f_{3}\right>$. Given this, we define a
ternary quadratic form $f_{\OO}$ 
associated to $\OO$ by
\[ f_{\OO}=d(\OO)\cdot N(X_{1}f_{1}+X_{2}f_{2}+X_{3}f_{3}). \]
Here multiplication by the (principal) ideal $d(\OO)$ is understood as
multiplication by a generator of $d(\OO)$. Hence, $f_{\OO}$ is only defined
up to multiplication by units in $R$.

This construction is due to Brzezinski, and it is a generalisation of
a construction originally made by Brandt in \cite{Brandt2}. This in
turn is a generalisation of a result in \cite{Latimer2}. The
construction of Brandt was investigated and clarified by Peters
in \cite{Peters}. Peters proved that it gives exactly all
Gorenstein orders under the restriction char$K\neq2$. This restriction
was eliminated  in \cite{Brzezinski7}. 
Notice that the modified correspondence given here  gives
all orders. Since there is no complete proof of the following theorem
in the literature, we will give one here which is almost self-contained. We
will closely follow the ideas in \cite{Brzezinski7}.

\thm{
Let $R$ be a principal ideal domain. The maps $f \mapsto C_{0}(f)$ and
$\OO \mapsto f_{\OO}$ are inverses to each 
other and the discriminants satisfy \mbox{$d(\OO)=(d(f_{\OO}))$.}
Furthermore, 
the maps give a bijection between similarity classes of 
non-degenerate ternary quadratic forms integral over $R$ and isomorphism
classes of quaternion $R$-orders. 
}

\begin{proof}
  First we prove that $C_{0}(f_{\OO})=\OO$. Suppose that
  $\OO=\left<x_{0},x_{1},x_{2},x_{3}\right>$ and
  $\OO^{\#}=\left<f_{0},f_{1},f_{2},f_{3}\right>$, where $\{f_{i}\}$ is
  the dual basis of $\{x_{i}\}$ and $x_{0}=1$. By definition, these
  bases satisfy 
  $Tr(x_{i}\bar{f_{j}})=\delta_{ij}$. In particular, $Tr(f_{0})=1$ and
  $Tr(f_{i})=0$, if $i>0$. It is straightforward to check that
  \begin{equation}
    \label{baseeq}
   x_{i}=Tr(f_{1}f_{2}f_{3})^{-1}\left( Tr(f_{j}f_{k}\bar{f_{0}})-
    f_{j}f_{k}\right),    
  \end{equation}
  where $(i,j,k)$ is an even permutation of $(1,2,3)$. For example,
  \begin{eqnarray*}
    Tr\left((Tr(f_{j}f_{k}\bar{f_{0}})-f_{j}f_{k})\bar{f_{i}}\right) & = & 
    Tr(f_{j}f_{k}\bar{f_{0}})Tr(\bar{f_{i}})-Tr(f_{j}f_{k}\bar{f_{i}})\\
     & = & Tr(f_{j}f_{k}f_{i})=Tr(f_{i}f_{j}f_{k})=Tr(f_{1}f_{2}f_{3}),
   \end{eqnarray*}
   since $(y_{1},y_{2},y_{3})\longmapsto Tr(y_{1}y_{2}y_{3})$ is
   trilinear and alternating on $\A_{0}$, and \mbox{$\bar{y}=-y$} on
   $\A_{0}$. By the proof of \cite[(3.2)]{Brzezinski7},
   $d=Tr(f_{1}f_{2}f_{3})^{-1}$ is a generator of $d(\OO)$ and
   $d f_{j}f_{k}$ is integral over $R$. Hence, we can rewrite
   \mref{baseeq} as
   \begin{equation}
     \label{baseeq2}
     df_{j}f_{k}=d\cdot Tr(f_{j}f_{k}\bar{f_{0}})-x_{i}
   \end{equation}
   and conclude that $d\cdot Tr(f_{j}f_{k}\bar{f_{0}})\in R$.

   From the definition of $f_{\OO}$, we get that
   \begin{equation} \label{formeq}
     f_{\OO}=d\cdot\sum_{i\leq j}a_{ij}X_{i}X_{j},
   \mbox{ where } a_{ij}=
   \left\{\begin{array}{ll}
     Tr(f_{i}\bar{f_{j}}), & \mbox{if } i\neq j \\
     N(f_{i}), & \mbox{if } i=j.
   \end{array} \right.
   \end{equation}
   Now it is straightforward to check that $\{e_{i}=df_{j}\bar{f_{k}}\}$,
   with $(i,j,k)$  an even permutation of $(1,2,3)$, satisfy the
   relations \mref{cliffrel} with $\{a_{ij}\}$ as in
   \mref{formeq}. Hence, we can identify $C_{0}(f_{\OO})$ with the
   order $\OO^{'}=\left<
   1,df_{1}\bar{f_{2}},df_{2}\bar{f_{3}},df_{3}\bar{f_{1}}
   \right>$. But since $df_{j}\bar{f_{k}}=-df_{j}f_{k}=x_{i}-
   d\cdot Tr(f_{j}f_{k}\bar{f_{0}})$, we get $\OO=\OO^{'}=C_{0}(f_{\OO})$.

   To prove the other direction, that is, $f_{C_{0}(f)}=f$ is a
   straightforward calculation. All one has to do is to determine the
   dual basis $f_{0},\,f_{1},\,f_{2},\,f_{3}$ of the basis satisfying
   \mref{cliffrel}, and then calculate
   $N(X_{1}f_{1}+X_{2}f_{2}+X_{3}f_{3})$. Of course, a computer with
   basic knowledge of non-commutative algebra is of great help. 
   
   The fact that isometric forms give isomorphic Clifford algebras
   follows directly from the universal property of Clifford
   algebras. Hence $f\sim g$ implies $C_{0}(f)\cong C_{0}(g)$, since
   as we remarked before, $C_{0}(f)=C_{0}(\epsilon f)$ if $\epsilon\in
   R^{*}$.

   Conversely, let $\OO_{1}$ and $\OO_{2}$ be isomorphic orders in
   $\A$. Then there 
   is an $a\in \A$ such that $\OO_{1}=a^{-1}\OO_{2}a$ and
   $\Lambda_{1}=a^{-1}\Lambda_{2}a$, where
   $\Lambda_{i}=\OO_{i}^{\#}\cap \A_{0}$. But $x\longmapsto a^{-1}xa$
   is an isometry of $\A_{0}$ with respect to $N_{0}$. Hence
   $f_{\OO_{1}}$ is isometric to $f_{\OO_{2}}$.

   Finally, it is an easy calculation to show that $d(\OO)=(d(f_{\OO}))$.
\end{proof}

\prop{
Let $f$ be a non-degenerate  ternary quadratic form integral over $R$ and
$\OO=C_{0}(f)$. If $f=b\cdot
g$, where $b\in R$ and $g$ is primitive, then the Brandt invariant of \OO\
is equal to $(b)$ and the Gorenstein closure of \OO\ is equal to
$C_{0}(g)$. In particular, $\OO$ is a Gorenstein order iff $f$ is
primitive. 
}

\begin{proof}
  It follows immediately from the relations \mref{cliffrel} that if
  $f=b\cdot g$, then $C_{0}(f)=R+b\cdot C_{0}(g)$. The proposition
  follows from this and \mref{gorenstein}.
\end{proof}

\section{Quaternion orders in the $\p$-adic case}\label{sec:loc}
In this section $K$ will be a $\p$-adic field with ring of integers
$R$ with prime $\p=(\pi)$. We will use the
correspondence in 
Section~\ref{sec:ter} and 
classifications of ternary quadratic forms to give a very explicit
classification of quaternion orders. It is much more elaborate in the
dyadic than in the non-dyadic case. Therefore we will restrict to
$2$-adic fields in the dyadic case, that is, fields in which $(2)$ is
prime. 

The discussion in this section will show that the maximal orders are
unique up to isomorphism in the $\p$-adic case. However, in the case of a
division algebra one can say even more. Therefore we include the
following well-known result \cite[ch. II]{Vigneras}:

\prop{ \label{maximal}
  Let $\OO$ be a maximal order in a quaternion algebra $\A$
over a local field $K$, and let $\p$ be  the maximal ideal
in $R$. If $\A \cong M(2,K)$, then $\OO$ is conjugate to $M(2,R)$,
so $d(\OO)=R$, and if $\A$ is a division algebra, 
then $\OO$ is unique and $d(\OO)=\p$.
}

We recall an invariant of orders in the local case introduced by
Eichler in \cite{Eichler2}. Let $k$ be the residue class field
of $K$ and $J(\OO)$ the Jacobson radical of $\OO$. If $\OO\not\cong
M_{2}(R)$, then the 
Eichler invariant $e(\OO)$ is defined to be:
\begin{equation}\label{eichlerdef}
e(\OO)= \left\{ \begin{array}{rl}
1, & \mbox{if } \OO / J(\OO) \cong k \oplus k \\
0,  & \mbox{if }  \OO / J(\OO) \cong k \\
-1,  & \mbox{if }  \OO / J(\OO) \mbox{ is a quadratic field extension of } k.
\end{array} \right. 
\end{equation}
Eichler also showed how to compute $e(\OO)$ easily. If
$x\in \A$, then the discriminant of 
$x$ is defined to be $\Delta(x)=Tr(x)^{2}-4N(x)$. The Eichler
invariant can be determined using the following result \cite[Satz
10]{Eichler2}:
\lemma{\label{calceinv}
The Eichler invariant $e(\OO)$ satisfies:
\begin{enumerate}
\item If $e(\OO)=0$, then $(\frac{\Delta(x)}{\p})=0 \;\forall x\in\OO$.
\item If $e(\OO)=1$, then $(\frac{\Delta(x)}{\p})\neq -1 \;\forall
  x\in\OO$, and $\exists x\in\OO : (\frac{\Delta(x)}{\p})=1$.
\item If $e(\OO)=-1$, then $(\frac{\Delta(x)}{\p})\neq 1 \;\forall
  x\in\OO$, and $\exists x\in\OO : (\frac{\Delta(x)}{\p})=-1$.
\end{enumerate}
Here $(\frac{a}{\p})=0$ iff $a\in\p$.
}
We remark that it is easy to show using \mref{gorenstein} and
\mref{calceinv} that if 
$\OO$ is not a Bass order, then $e(\OO)=0$. 

Let $f$ be a ternary quadratic form integral over $R$. We can restrict
to $f$ primitive, since if $f_{1}=b_{1}\cdot g_{1}$ and
$f_{2}=b_{2}\cdot g_{2}$ with $b_{i}\in R$ and $g_{i}$ primitive
forms, then
$f_{1}\sim f_{2}$ iff $g_{1}\sim g_{2}$ and the ideals $(b_{1})$ and
$(b_{2})$ are equal. There is a classification of quadratic lattices
in the local case in \cite[\S\S 92,93]{O'Meara}, which we will make use
of below. 

In the non-dyadic case every form is isometric to a diagonal form, but
this is not true in the $2$-adic case. We define the $2$-dimensional
quadratic forms $H$ and $J$ with the corresponding matrices
\[ H=\Matrix{0}{1}{1}{0} \mbox{ and } J=\Matrix{2}{1}{1}{2}. \]
Since we may multiply $f$ by elements in $R^{*}$,  we can conclude the
following from the results in \cite{O'Meara}:

\prop{ \label{standard}
Let $f$ be a primitive ternary quadratic form over $R$.\\ If $\p$ is a
non-dyadic prime, then 
  \[ f\sim \left<1\right>\perp\left<\delta\pi^{r}\right>\perp\left<
    \epsilon\pi^{s}\right>. \]
If $\p=(2)$, then $f\sim f_{i}$ for some $1\leq i\leq5$, where
  \begin{eqnarray*}
    f_{1} & = & \left<1\right>\perp\left<\delta 2^{r}\right>\perp\left< 
    \epsilon 2^{s}\right>, \\
    f_{2}(r) & = & \left<1\right>\perp 2^{r}H, \\
    f_{3}(r) & = & H\perp \left<2^{r}\right>, \\
    f_{4}(r) & = & \left<1\right>\perp 2^{r}J, \\
    f_{5}(r) & = & J\perp \left<2^{r}\right>. 
  \end{eqnarray*}
Here $\epsilon,\delta\in R^{*}$ and $0\leq r\leq s$.
}

The quadratic forms in \mref{standard} will be called standard
forms. 
\remn{
The only cases when any of the non-diagonal forms $f_{2},\ldots,f_{5}$ are
similar to a diagonal form are 
\begin{equation} \label{nondia1}
\left<1\right>\perp 2H\sim\ter11{-1} \mbox{ and } 
   \left<1\right>\perp 2J\sim \ter111.
\end{equation}
Furthermore, we also have
\begin{equation}
  \label{nondia2}
  f_{2}(0)=f_{3}(0)\sim f_{4}(0)=f_{5}(0) \mbox{ and } f_{2}(2)\sim
  f_{4}(2).
\end{equation}
However, $f_{i}\not\sim f_{j}$ if $i\neq j$ in all cases except
\mref{nondia1} and \mref{nondia2}. 
}

The following proposition gives simple criteria on $f$ whether
$C_{0}(f)$ is a Bass order or not.
\prop{ \label{basscrit}
If $\p$ is non-dyadic and $f\sim
\left<1\right>\perp\left<\delta\pi^{r}\right>\perp\left< 
  \epsilon\pi^{s}\right>$ with $r\leq s$, then $C_{0}(f)$ is a Bass
order iff $r\leq 1$.

If $\p=(2)$, then
$C_{0}(f)$ is a Bass order iff $f$ is similar to any of the forms 
\[ \begin{array}{l}
  \ter{1}{\delta}{\epsilon2^{r}},\mbox{ with } \delta\md14  \mbox{ or
  } r\leq1, \\
  \ter{1}{\delta2}{\epsilon2^{r}}, \\
  f_{3}(r) \mbox{ or }   f_{5}(r).
\end{array} \]
}
\begin{proof}
  It follows from \mref{norm} that if
  $f=\ter{1}{\delta\pi^{r}}{\epsilon\pi^{s}}$, then the norm form in
  $\OO=C_{0}(f)$ is given by
  \[ N=\left<1\right>\perp
  \ter{\delta\pi^{r}}{\epsilon\pi^{s}}{\delta\epsilon\pi^{r+s}}. \]
  It is easy to check from this that $\OO$ contains a primitive element (and
  hence is a Bass order) exactly in the cases in the proposition. The
  non-diagonal cases are analogous. 
\end{proof}

If we have a standard form $f$, then it is easy to calculate the
Eichler invariant of $C_{0}(f)$. A direct computation using
\mref{calceinv} gives:
\prop{
If $\OO=C_{0}(f)$, then the Eichler invariant $e(\OO)$ satisfies:

If $\p$ is non-dyadic and $f\sim
    \left<1\right>\perp\left<\delta\pi^{r}\right>\perp\left< 
    \epsilon\pi^{s}\right>$ with $r\leq s$, then:
  \begin{eqnarray*}
  e(\OO)=1 & \mbox{iff} & r=0,\,s\geq 1 \mbox{ and }
  (\frac{-\delta}{\p})=1. \\
  e(\OO)=-1 & \mbox{iff} & r=0,\,s\geq 1 \mbox{ and }
  (\frac{-\delta}{\p})=-1.
  \end{eqnarray*}

If $\p=(2)$, then:
  \begin{eqnarray*}
  e(\OO)=1 & \mbox{iff} & f\sim H\perp \left<2^{r}\right>,\, r\geq 1. \\
  e(\OO)=-1 & \mbox{iff} & f\sim J\perp \left<2^{r}\right>,\, r\geq 1.
  \end{eqnarray*}
}

There are descriptions of all orders in quaternion algebras over local
fields in both \cite{Brzezinski2} and \cite{Korner1}. In
\cite{Brzezinski2}, which is more ring-theoretic in nature, there are
also detailed information on relations between orders. Now we will, as
a complement to these descriptions, give a set of primitive ternary
quadratic forms $M_{1}$ such that $M^{'}=\{C_{0}(f) : f\in M_{1}\}$
is a set of representatives of all isomorphism classes of Gorenstein
orders. Then 
\[ M=\{C_{0}(\pi^{b}f):f\in M_{1},\,b\geq0 \}\]
 is a set of
representatives of all isomorphism classes of orders. One of the
advantages of our description is that it is very well suited for
explicit calculations, since $f$ in standard form gives a convenient
basis of $C_{0}(f)$.

We start with the non-dyadic case. If
\[ f_{1}=\ter{1}{\delta_{1}\pi^{r_{1}}}{\epsilon_{1}\pi^{s_{1}}} \sim
f_{2}=\ter{1}{\delta_{2}\pi^{r_{2}}}{\epsilon_{2}\pi^{s_{2}}}, \]
then $r_{1}=r_{2}$ and $s_{1}=s_{2}$. Conversely, if $r_{1}=r_{2}$,
$s_{1}=s_{2},\,(\frac{\delta_{1}}{\p})=(\frac{\delta_{2}}{\p})$ and
$(\frac{\epsilon_{1}}{\p})=(\frac{\epsilon_{2}}{\p})$, then $f_{1}\sim
f_{2}$. Hence,  we have at most $4$ different
classes given $r$ and $s$. But if $r=0$ or $r=s$, then the number of
classes is divided 
by $2$ \cite[92:1]{O'Meara}. This is summarised in:
\prop{ \label{nonclass}
Let $\p$ be a non-dyadic prime and let $\OO$ be a Gorenstein order in a
quaternion algebra over $K$. Then $\OO\cong C_{0}(f)$, where $f$ is
uniquely chosen among the following quadratic forms:
\[ \begin{array}{llll}
1. & \left<1\right>\perp\left<1\right>\perp\left<1\right> &
3. & \left<1\right>\perp\left<\pi^{r}\right>\perp
\left<\epsilon_{1}\pi^{r}\right>,\; r\geq 1 \\
2. & \left<1\right>\perp\left<-\epsilon_{1}\right>\perp
\left<\pi^{s}\right>,\; s\geq 1 &
4. & \left<1\right>\perp\left<\epsilon_{1}\pi^{r}\right>\perp
\left<\epsilon_{2}\pi^{s}\right>,\; s>r\geq 1 
\end{array} \]
Here $\epsilon_{1}$ and $\epsilon_{2}$ are to be chosen arbitrarily in
$\{1,\delta\}$, where $(\frac{\delta}{\p})=-1$. 
}

With \mref{basscrit} and \mref{nonclass} at hand, we are able to
determine the number of isomorphism classes of Bass and Gorenstein
orders with given discriminant. Furthermore, if $t(n)$ is the number
of isomorphism classes of orders with discriminant $\p^{n}$ and $g(n)$
is the number of isomorphism classes of Gorenstein orders with
discriminant $\p^{n}$, then
\begin{equation}
  \label{linrec}
  t(n)=g(n)+t(n-3). 
\end{equation}
The reason for this is that an order \OO\  is either Gorenstein or else
$\OO=R+\p\OO^{'}$ with $[\OO^{'}:\OO]=\p^{3}$. Moreover, if
$\OO_{1}=R+\p\OO^{'}_{1}$ and $\OO_{2}=R+\p\OO^{'}_{2}$, then
$\OO_{1}\cong\OO_{2}$ iff $\OO_{1}^{'}\cong\OO_{2}^{'}$. The linear
recursion \mref{linrec} is easy to solve, since $g(n)$ proves to be
more or less a 
linear function. We summarise everything in  Table~\ref{tab:non}, in
which we separate orders in $\h_{\p}$ and $M_{2}=M_{2}(K)$. The two
variables $n$ and $\alpha$ in  Table~\ref{tab:non} satisfy
\[ n\geq3 \mbox{ and } \alpha=\left\{
\begin{array}{ll}
  \frac13, & \mbox{if } n\md03 \\
  0, & \mbox{otherwise.}
\end{array}\right.
\]

\begin{table}[hbtp]  
\[
\begin{array}{|c|cc||cc||cc|}\hline
   & \multicolumn{2}{c||}{\mbox{BASS}} &
 \multicolumn{2}{c||}{\mbox{GOR.}}  & \multicolumn{2}{c|}{\mbox{TOTALLY}}
 \\ \hline 
 \mbox{Disc.} & \h_{\p} & M_{2} & \h_{\p} & M_{2}
 & \h_{\p} & M_{2} \\ \hline
 1 & 0 & 1 & 0 & 1 & 0 & 1  \\
 \p & 1 & 1 & 1 & 1 & 1 & 1  \\
 \p^{2} & 1 & 3 &  1 & 3 & 1 & 3 \\
 \p^{n},n \mbox{ odd} & 3 & 3 & n & n 
 &\frac18(n^{2}+4n+3) & \frac{1}{24}(5n^{2}+12n+7)+\alpha\\
 \p^{n},n \mbox{ even} & 2 & 4 & \frac{n}{2} & \frac{3n}{2} 
 &\frac18(n^{2}+2n) & \frac{1}{24}(5n^{2}+18n+16)+\alpha \\ \hline 
\end{array}
\]
  \caption{The number of isomorphism classes of Bass, Gorenstein and
  arbitrary orders in 
  $\h_{\p}$ and $M_2=M_2(K)$, when $\p$ is non-dyadic.} 
  \label{tab:non}
\end{table}

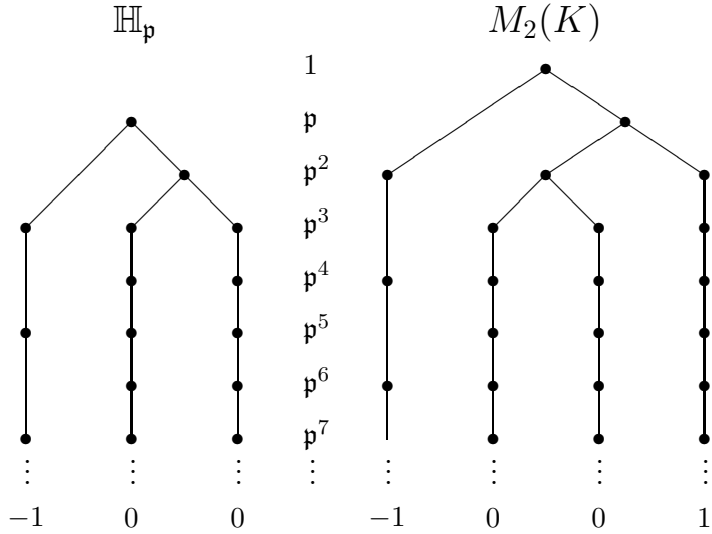
\begin{figure}[hbtp]
  \begin{picture}(320,215)(-40,-20)

\multiput(20,20)(0,40){3}{\circle*{4}}
\multiput(60,20)(0,20){5}{\circle*{4}}
\multiput(100,20)(0,20){5}{\circle*{4}}
\put(60,140){\circle*{4}}
\put(80,120){\circle*{4}}
\multiput(20,20)(40,0){3}{\line(0,1){80}}
\put(20,100){\line(1,1){40}}
\put(60,100){\line(1,1){20}}
\put(100,100){\line(-1,1){40}}

\put(19,3){$\vdots$}
\put(59,3){$\vdots$}
\put(99,3){$\vdots$}
\put(127,3){$\vdots$}
\put(20,-10){\makebox(0,0){$-1$}}
\put(60,-10){\makebox(0,0){$0$}}
\put(100,-10){\makebox(0,0){$0$}}

\put(125,18){$\p^{7}$} 
\put(125,38){$\p^{6}$} 
\put(125,58){$\p^{5}$} 
\put(125,78){$\p^{4}$} 
\put(125,98){$\p^{3}$} 
\put(125,118){$\p^{2}$} 
\put(125,138){$\p$} 
\put(125,158){$1$} 
\put(53,174){\Large $\h_{\p}$}

\put(133,20){
  \begin{picture}(150,300)
\multiput(20,20)(0,40){3}{\circle*{4}}
\multiput(60,0)(0,20){5}{\circle*{4}}
\multiput(100,0)(0,20){5}{\circle*{4}}
\multiput(140,0)(0,20){6}{\circle*{4}}
\put(80,140){\circle*{4}}
\put(80,100){\circle*{4}}
\put(110,120){\circle*{4}}
\multiput(20,0)(120,0){2}{\line(0,1){100}}
\multiput(60,0)(40,0){2}{\line(0,1){80}}
\put(20,100){\line(3,2){60}}
\put(60,80){\line(1,1){20}}
\put(100,80){\line(-1,1){20}}
\put(80,100){\line(3,2){30}}
\put(140,100){\line(-3,2){30}}
\put(110,120){\line(-3,2){30}}

\put(19,-17){$\vdots$}
\put(59,-17){$\vdots$}
\put(99,-17){$\vdots$}
\put(139,-17){$\vdots$}
\put(20,-30){\makebox(0,0){$-1$}}
\put(60,-30){\makebox(0,0){$0$}}
\put(100,-30){\makebox(0,0){$0$}}
\put(140,-30){\makebox(0,0){$1$}}

\put(58,154){\Large $M_{2}(K)$}    
  \end{picture}
}

\end{picture}
  \begin{center}
  \caption{The trees of isomorphism classes of Bass orders in the
    non-dyadic case.}  \label{fig:non}
  \end{center}  
\end{figure}

\vspace{-3mm}
If we restrict to Bass orders, then we can also easily draw
conclusions on the relations between the classes of orders by
using results first proved by Eichler in \cite[Satz
12]{Eichler2}. These relations on inclusions between orders were
generalised in 
\cite{Brzezinski2}. We illustrate this by the trees in
Figure~\ref{fig:non}. Every 
node in a tree represents an isomorphism class of Bass orders.
Different isomorphism classes with the same discriminant are on the
same level. There 
is an edge between nodes $n_{1}$ and $n_{2}$ iff given $\OO_{1}\in
n_{1}$, then there exists $\OO_{2}\in n_{2}$ such that $\OO_{1}$ is a
maximal suborder in $\OO_{2}$ or vice versa. The numbers at the bottom
are the Eichler invariants of the orders in that column. Moreover, the
order \OO\ in $\h_{\p}$ with $d(\OO)=\p$ is the maximal order and has
$e(\OO)=-1$, and the orders in $M_{2}(K)$ with $d(\OO)=\p$ have $e(\OO)=1$.

We now turn to the $2$-adic case. The only relations for the
non-diagonal forms are given by \mref{nondia1} and \mref{nondia2}. For
the diagonal forms, we have as in the non-dyadic case that 
\[ f_{1}=\ter{1}{\delta_{1}\pi^{r_{1}}}{\epsilon_{1}\pi^{s_{1}}} \sim
f_{2}=\ter{1}{\delta_{2}\pi^{r_{2}}}{\epsilon_{2}\pi^{s_{2}}}, \]
implies that $r_{1}=r_{2}$ and $s_{1}=s_{2}$. Conversely, if $r_{1}=r_{2}$,
$s_{1}=s_{2}$, \mbox{$\delta_{1}\md{\delta_{2}}8$} and
$\epsilon_{1}\md{\epsilon_{2}}8$, then $f_{1}\sim
f_{2}$. Hence, we have at most $16$ classes for every pair
$(r,s)$. But if $r<3$ or $s-r<3$, then the number of classes is
reduced. We have used the `canonical form' in \cite{BJones} rather
than \cite{O'Meara} in order to get the following set of representatives:

\prop{ \label{dyadicclass}
Let $\p=(2)$ and let $\OO$ be a Gorenstein order in a
quaternion algebra over $K$. Then $\OO\cong C_{0}(f)$, where $f$ is
uniquely chosen among the  quadratic forms in Table~$\ref{dyadicforms}$.
}

\begin{table}[hbtp]{\small 
\[ \begin{array}{|rlrl|}\hline
&&&\\
1. & H\perp \left<2^{r}\right>,\, r\geq 0 &
12. & \ter{1}{\delta_{4}2^{2}}{\delta_{4}2^{r}},\, r\geq 2 \\
2. & J\perp \left<2^{r}\right>,\, r\geq 1 &
13. & \ter{1}{2^{2}}{3\cdot2^{r}},\, r\geq 2 \\
3. & \left<1\right>\perp 2^{r}H,\, r\geq 2 &
14. & \ter{1}{3\cdot2^{2}}{2^{r}},\, r\geq 3 \\
4. & \left<1\right>\perp 2^{r}J,\, r\geq 3 &
15. & \ter{1}{\delta_{7}2^{2}}{7\cdot2^{r}},\, r\geq 5 \\
5. & \ter{1}{1}{\delta_{1}2^{r}},\, r\geq 0 &
16. & \ter{1}{\delta_{2}2^{2}}{2^{r}},\, r\geq 5 \\
6. & \ter{1}{3}{2^{r}},\, r\geq 2 &
17. & \ter{1}{2^{r}}{\delta_{5}2^{r+s}},\, r\geq 3,\,s\geq0  \\
7. & \ter{1}{7}{2^{r}},\, r\geq 3 &
18. & \ter{1}{7\cdot2^{r}}{\delta_{6}2^{r+s}},\, r\geq 3,\,s\geq0  \\
8. & \ter{1}{5}{\delta_{1}2^{r}},\, r\geq 3 &
19. & \ter{1}{3\cdot2^{r}}{\delta_{3}2^{r+s}},\, r\geq 3,\,s\geq1 \\
9. & \ter{1}{3\cdot2}{\delta_{1}2^{r}},\, r\geq 1 &
20. & \ter{1}{7\cdot2^{r}}{\delta_{3}2^{r+s}},\, r\geq 3,\,s\geq3  \\
10. & \ter{1}{2}{\delta_{3}2^{r}},\, r\geq 3 &
21. & \ter{1}{3\cdot2^{r}}{\delta_{6}2^{r+s}},\, r\geq 3,\,s\geq3 \\
11. & \ter{1}{\delta_{2}2}{\delta_{3}2^{r}},\, r\geq 4  &
22. & \ter{1}{5\cdot2^{r}}{\delta_{5}2^{r+s}},\, r\geq 3,\,s\geq3 \\ 
&&&\\ \hline
\end{array} \]
}
  \caption{Representatives of the similarity classes of primitive
    ternary quadratic forms over  $2$-adic integers. The
    $\delta_{i}$:s are to be chosen arbitrarily 
    according to $\delta_{1}\in\{1,3\},\,\delta_{2}\in\{5,7\},\, 
    \delta_{3}\in\{1,5\}$,
    $\delta_{4}\in\{1,7\},\,\delta_{5}\in\{1,3,5,7\} 
    ,\,\delta_{6}\in\{3,7\}$ and $\delta_{7}\in\{3,5\}$.} 
  \label{dyadicforms}
\end{table}

For example 11. in Table~\ref{dyadicforms} gives rise to $4$ classes
  for every $r\geq4$, namely 
  $(\delta_{2},\delta_{3})\in\{(5,1),\,(5,5),\,(7,1),\,(7,5)\}$. We get
  from \mref{basscrit}, that the Bass orders are those in 
\{1,\,2,\,5,\,8,\,9,\,10,\,11\}. Now with \mref{dyadicclass} at our
  disposal, we 
can determine the number of isomorphism classes of orders also in the
$2$-adic case. We summarise everything in  Table~\ref{tab:dya}. The
  numbers $n$ and 
$\beta$ in  Table~\ref{tab:dya} satisfy
\[ n\geq 9 \mbox{ and } \beta=\left\{
  \begin{array}{ll}
    \frac13, & \mbox{if } n\md23 \\
    0, & \mbox{otherwise.}
  \end{array} \right.
\]

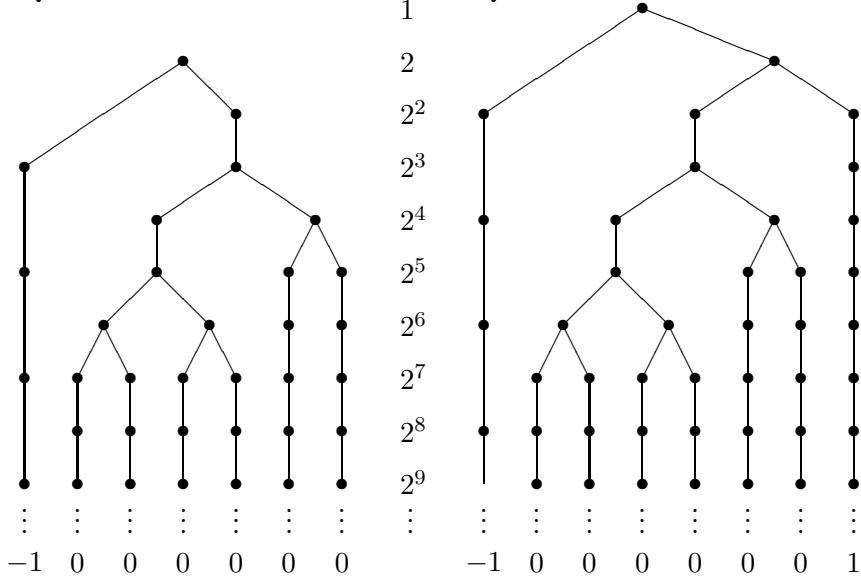
\begin{figure}[!hbtp]
  \begin{picture}(320,200)(-10,-10)

\multiput(20,20)(0,40){4}{\circle*{4}}
\multiput(40,20)(0,20){3}{\circle*{4}}
\multiput(60,20)(0,20){3}{\circle*{4}}
\multiput(80,20)(0,20){3}{\circle*{4}}
\multiput(100,20)(0,20){3}{\circle*{4}}
\multiput(120,20)(0,20){5}{\circle*{4}}
\multiput(140,20)(0,20){5}{\circle*{4}}
\put(50,80){\circle*{4}}
\put(90,80){\circle*{4}}
\put(70,100){\circle*{4}}
\put(70,120){\circle*{4}}
\put(130,120){\circle*{4}}
\put(100,140){\circle*{4}}
\put(100,160){\circle*{4}}
\put(80,180){\circle*{4}}
\multiput(40,20)(20,0){4}{\line(0,1){40}}
\multiput(120,20)(20,0){2}{\line(0,1){80}}
\put(20,20){\line(0,1){120}}
\put(20,140){\line(3,2){60}}
\put(40,60){\line(1,2){10}}
\put(60,60){\line(-1,2){10}}
\put(80,60){\line(1,2){10}}
\put(100,60){\line(-1,2){10}}
\put(50,80){\line(1,1){20}}
\put(90,80){\line(-1,1){20}}
\put(70,100){\line(0,1){20}}
\put(70,120){\line(3,2){30}}
\put(130,120){\line(-3,2){30}}
\put(120,100){\line(1,2){10}}
\put(140,100){\line(-1,2){10}}
\put(100,140){\line(0,1){20}}
\put(100,160){\line(-1,1){20}}

\multiput(20,9)(20,0){7}{\makebox(0,0){$\vdots$}}
\put(166,9){\makebox(0,0){$\vdots$}}
\put(20,-10){\makebox(0,0){$-1$}}
\multiput(40,-10)(20,0){6}{\makebox(0,0){$0$}}

\put(162,16){$2^{9}$} 
\put(162,36){$2^{8}$} 
\put(162,56){$2^{7}$} 
\put(162,76){$2^{6}$} 
\put(162,96){$2^{5}$} 
\put(162,116){$2^{4}$} 
\put(162,136){$2^{3}$} 
\put(162,156){$2^{2}$} 
\put(162,176){$2$} 
\put(162,196){$1$} 
\put(80,212){\makebox(0,0){\Large $\h_{2}$}}

\put(170,20){
  \begin{picture}(240,300)
\multiput(20,20)(0,40){4}{\circle*{4}}
\multiput(40,0)(0,20){3}{\circle*{4}}
\multiput(60,0)(0,20){3}{\circle*{4}}
\multiput(80,0)(0,20){3}{\circle*{4}}
\multiput(100,0)(0,20){3}{\circle*{4}}
\multiput(120,0)(0,20){5}{\circle*{4}}
\multiput(140,0)(0,20){5}{\circle*{4}}
\multiput(160,0)(0,20){8}{\circle*{4}}
\put(50,60){\circle*{4}}
\put(90,60){\circle*{4}}
\put(70,80){\circle*{4}}
\put(70,100){\circle*{4}}
\put(130,100){\circle*{4}}
\put(100,120){\circle*{4}}
\put(100,140){\circle*{4}}
\put(130,160){\circle*{4}}
\put(80,180){\circle*{4}}
\multiput(40,0)(20,0){4}{\line(0,1){40}}
\multiput(120,0)(20,0){2}{\line(0,1){80}}
\put(20,0){\line(0,1){140}}
\put(160,0){\line(0,1){140}}
\put(40,40){\line(1,2){10}}
\put(60,40){\line(-1,2){10}}
\put(80,40){\line(1,2){10}}
\put(100,40){\line(-1,2){10}}
\put(50,60){\line(1,1){20}}
\put(90,60){\line(-1,1){20}}
\put(70,80){\line(0,1){20}}
\put(120,80){\line(1,2){10}}
\put(140,80){\line(-1,2){10}}
\put(70,100){\line(3,2){30}}
\put(130,100){\line(-3,2){30}}
\put(100,120){\line(0,1){20}}
\put(100,140){\line(3,2){30}}
\put(160,140){\line(-3,2){30}}
\put(130,160){\line(-5,2){50}}
\put(20,140){\line(3,2){60}}

\multiput(20,-11)(20,0){8}{\makebox(0,0){$\vdots$}}
\put(20,-30){\makebox(0,0){$-1$}}
\multiput(40,-30)(20,0){6}{\makebox(0,0){$0$}}
\put(160,-30){\makebox(0,0){$1$}}

\put(80,192){\makebox(0,0){\Large $M_{2}(K)$}}   
  \end{picture}
}

\setlength{\unitlength}{1pt}

\end{picture}
  \begin{center}
    \caption{The trees of isomorphism classes of Bass orders in the
    $2$-adic case.}\label{fig:dya}
  \end{center}
\end{figure}

The trees of isomorphism classes of Bass orders in the $2$-adic case
are shown in Figure~\ref{fig:dya}.

\begin{center}
\begin{table}[hbtp]  
{\scriptsize
\[
\begin{array}{|c|cc||cc||cc|}\hline
   & \multicolumn{2}{c||}{\mbox{BASS}} &
 \multicolumn{2}{c||}{\mbox{GORENSTEIN}}  &
 \multicolumn{2}{c|}{\mbox{TOTALLY}} 
 \\ \hline 
 \mbox{Disc.} & \h_{2} & M_{2} & \h_{2} & M_{2}
 & \h_{2} & M_{2} \\ \hline
 1 & 0 & 1 & 0 & 1 & 0 & 1  \\
 2 & 1 & 1 & 1 & 1 & 1 & 1  \\
 2^{2} & 1 & 3 &  1 & 3 & 1 & 3 \\
 2^{3} & 2 & 2 &  2 & 2 & 2 & 3 \\
 2^{4} & 2 & 4 &  2 & 6 & 3 & 7 \\
 2^{5} & 4 & 4 &  5 & 5 & 6 & 8 \\
 2^{6} & 4 & 6 &  6 & 11 & 8 & 14 \\
 2^{7} & 7 & 7 &  10 & 10 & 13 & 17 \\
 2^{8} & 6 & 8 & 10 & 18 & 16 & 26 \\
 2^{n},n \mbox{ odd} & 7 & 7 &  4(n-5) & 4(n-5) &
 \frac{1}{12}(7n^{2}-46n+135)+\beta &
 \frac14(3n^{2}-22n+75)+3\beta \\
 2^{n},n \mbox{ even} & 6 & 8 & 3(n-5) & 5(n-5) & 
 \frac{1}{12}(7n^{2}-52n+156)+\beta &
 \frac14(3n^{2}-20n+68)+3\beta \\ \hline 
\end{array}
\]
}
  \caption{The number of isomorphism classes of Bass, Gorenstein and
  arbitrary orders in 
  $\h_{\p}$ and $M_2=M_2(K)$, when $\p=(2)$.} 
  \label{tab:dya}
\end{table}
\end{center}

\section{Quaternion orders in the number field case}\label{sec:glo}
In this section, $K$ will be an algebraic number field with ring of
integers $R$. 

A powerful tool when studying orders in quaternion algebras over
algebraic number 
fields is to consider completions $\OO_{\p}=\OO\otimes_{R}R_{\p}$ with
respect to the  primes ideals in $R$. The basic result which makes
this so useful is the following 
local-global correspondence \cite[prop. III.5.1]{Vigneras}:

\prop{ \label{localglobal}
  If $\Lambda^{'}$ is a lattice on $\A$, then there is a bijection 
between $\{$lattices on $\A \}$ and $\{ (\Lambda_{\p})_{\p \in \Omega_{f}}:
\Lambda_{\p}$ lattice on $\A_{\p}$, $\Lambda_{\p}=\Lambda^{'}_{\p}$ for almost
all $\p \}$ given by 
\[ \Lambda \mapsto (\Lambda_{\p}) \mbox{ and } (\Lambda_{\p}) \mapsto
\Lambda= \{x \in \A: 
x \in \Lambda_{\p}, \: \forall \p\in \Omega_{f} \}. \]
}

The local-global principle remains true if we restrict to orders, since
$\Lambda$ is an order iff $\Lambda_{\p}$ is an order for all
$\p\in\Omega_{f}$. It is immediate from the
definitions that the discriminant and index satisfy
\begin{equation}
  \label{lgdisc}
  d(\OO)_{\p}=d(\OO_{\p}) \mbox{ and } \left[\OO^{'}_{\p} :
  \OO_{\p}\right]=\left[\OO^{'} : \OO\right]_{\p}. 
\end{equation}

A direct consequence of \mref{localglobal} is that \OO\ is maximal iff
$\OO_{\p}$ is maximal for all $\p\in\Omega_{f}$. From this,
\mref{maximal} and \mref{lgdisc}, we get that 
\begin{equation}
  \label{maxcond}
  \OO \mbox{ is maximal in }\A \mbox{ iff } d(\OO)=d(\A).
\end{equation}
It is also clear that $\OO$ is a Bass (Gorenstein) order iff $\OO_{\p}$
is a Bass (Gorenstein) order for all $\p\in\Omega_{f}$.

As already mentioned, the local-global principle does not apply to
isomorphism classes since $\OO_{\p}\cong\OO^{'}_{\p}\;
\forall\p\in\Omega_{f}$ does not imply $\OO\cong\OO^{'}$ in general. We
say that \OO\ and $\OO^{'}$ are in the same genus if $\OO_{\p}\cong\OO^{'}_{\p}\;
\forall\p\in\Omega_{f}$, and that they are of the same type if
$\OO\cong\OO^{'}$. The number of types (isomorphism classes) in the genus of
$\OO$ is called the type number of \OO\ and will be denoted by $t(\OO)$.

Closely related to $t(\OO)$ are the one- and two-sided class numbers of
\OO, which are defined as follows. Let $L(\OO)$ be the set of locally
principal left \OO-ideals. We define an equivalence relation on $L(\OO)$
by $\Lambda_{1}\sim\Lambda_{2}$ iff there exists $a\in\A$ such that
$\Lambda_{1}=\Lambda_{2}a$. The number of equivalence classes in
$L(\OO)$ with respect to this relation is the one-sided class number
$h(\OO)$ of \OO. In the literature this is often simply referred to as
the  class number of \OO. The facts that $h(\OO)<\infty$ and that we get
the same number of classes if we take right ideals
instead were first proved by Brandt in \cite{Brandt1}. That paper is
the first general investigation of the ideal theory of quaternion
algebras.  

If we restrict to two-sided ideals instead
and do the same construction, then the number of equivalence classes
is called the two-sided class number of \OO\ and will be denoted by
$H(\OO)$. 

\remn{
There exist orders with ideals which are not locally principal. For
example, let $\p=(\pi)$ be a principal ideal in $R$ and
\[ \OO=(\pi,\pi)_{R}. \] 
If $\Lambda=\{x\in\OO : N(x)\in\p\}=\p+Ri+Rj+Rij$, then clearly
$\Lambda$ is a two-sided \OO-ideal. Furthermore $[\OO : \Lambda]=\p$ and
$\Lambda=\OO i+\OO j$. 

If $g=\pi\alpha+\beta i+\gamma j+\delta ij\in\Lambda_{\p}$, then
\[ \left( 
  \begin{array}{c}
    g\cdot 1 \\
    g\cdot i \\
    g\cdot j \\
    g\cdot ij
  \end{array} \right)
=\left( \begin{array}{cccc} \pi\alpha & \beta & \gamma & \delta \\
    \pi\beta & \pi\alpha & -\delta & -\gamma \\
    \pi\gamma & \pi\delta & \pi\alpha & \beta  \\
    -\pi^{2}\delta &-\pi\gamma & \pi\beta & \pi\alpha
  \end{array} \right)
  \left( 
  \begin{array}{c}
     1 \\
     i \\
     j \\
     ij
  \end{array} \right)=:A
  \left( 
  \begin{array}{c}
     1 \\
     i \\
     j \\
     ij
   \end{array} \right).
 \]
Hence, $[\OO_{\p} : \OO_{\p} g]=(\det(A))\subseteq\p^{2}$ and this implies that
$\Lambda_{\p}$ is not principal since $[\OO_{\p} : \Lambda_{\p}]=\p$.
}

Let $\Lambda_{1},\ldots,\Lambda_{h(\OO)}$ be a set of representatives of
left \OO-ideal classes. Then every order in the genus of
\OO\ is isomorphic to at least one of the right orders
$\OO_{r}(\Lambda_{i})$. In fact, if $\OO^{'}$ is in the genus of \OO,
then $\OO^{'}$ is isomorphic to exactly $H(\OO^{'})$ of the orders
$\OO_{r}(\Lambda_{i})$. If we fill in the details, we get a proof of
the following proposition. For an ad\'elic proof, see
\cite[p. 88]{Vigneras}. 
\prop{
Let \OO\ be an arbitrary order in a quaternion algebra over $K$. Then
\[ h(\OO)=\sum_{i=1}^{t(\OO)}H(\OO_{i}), \]
where $\OO_{1},\ldots,\OO_{t(\OO)}$ are a set of representatives of the
types in the genus of \OO. In particular, $h(\OO)=h(\OO^{'})$, if $\OO$ and
$\OO^{'}$ are in the same genus.
}

For  results regarding the calculation of class and
type numbers, see \cite{Brzezinski8}, \cite{Brzezinski1},
\cite{Eichler4}, \cite{Eichler1}, \cite{Korner2}, \cite{Pizer2},
\cite{Pizer3} and [SJ2]. 

\remn{
There is a little confusion on notions of `class number' regarding
quadratic forms versus quaternion orders. Let $N$ be the norm form of
an order $\OO$. The class number of the quadratic form $N$ is
defined to be the number of isometry classes in the genus of $N$,
that is, forms locally isometric to $N$. Hence, the class number of $N$
is equal to the type number of $\OO$.
}

\exn{
In order to illustrate our discussion, we will as an example determine the
number of isomorphism classes of orders \OO\ in rational quaternion
algebras with  $d(\OO)=(72)=(2^{3}\cdot3^{2})$. If $\OO$
is an order 
in an algebra $\A$, then $d(\A)$ divides $d(\OO)$. Since the discriminant
determines the rational quaternion algebra and every positive square
free integer is the discriminant of a quaternion algebra, we get $4$
possibilities: $d(\A)\in\{1,2,3,6\}$. The cases $d(\A)\in\{2,3\}$ are
definite algebras and the other two are indefinite. 

From
Table~\ref{tab:non}, we get that there are $1$ isomorphism class of
orders with $d(\OO_{3})=(72)=(3^{2})$ in $\h_{3}$ and $3$ in
$M_{2}(\q_{3})$. Furthermore, Table~\ref{tab:dya} implies that there
are $2$ isomorphism classes of orders with $d(\OO_{2})=(72)=(2^{3})$ in
$\h_{2}$ and $3$ in $M_{2}(\q_{2})$. Of all these orders only one in
$M_{2}(\q_{2})$ is
not Gorenstein or Bass. Now from \mref{localglobal}, we derive the
second column of Table~\ref{tab:ex1}. 
\begin{table}[hbtp]
  \[ \begin{array}{|c|c|c|} \hline 
    d(\A) & \mbox{\#genera} & \mbox{\#types} \\ \hline
    1 & 9 & 9 \\
    2 & 6 & 10 \\
    3 & 3 & 3 \\
    6 & 2 & 2 \\ \hline
  \end{array} \]
    \caption{The number of genera and types of rational orders with
    discriminant equal to (72).} 
    \label{tab:ex1}
\end{table}

It is easy to show, using for
example the results in \cite{Brzezinski1}, that the type numbers of the
orders in the two indefinite algebras are all equal to $1$ in our
case. The easiest way to determine  the type numbers of the definite
orders is to check the tables in \cite{BI}. These tables reveal that
$4$ of the genera in the algebra ramified at $2$ have $2$ classes and
all other genera only one.
This gives the last column of Table~\ref{tab:ex1}. Of course, the
tables in \cite{BI} give all this information concerning the definite
algebras in this case. 
}

We conclude by giving an explicit basis of a maximal
order in an arbitrary quaternion algebra \A\ with $d(\A)=(d)$
principal. First choose a generator $a$ of a prime ideal satisfying
\mref{alpha}, so that $\A\cong(a,-d)_{K}$. We start with the order
$\OO^{'}=(a,-d)_{R}$, which has $d(\OO^{'})=(4ad)$. From \mref{discindex} and
\mref{maxcond}, we get that an order $\OO\supseteq\OO^{'}$ is maximal iff
$[\OO : \OO^{'}]=(4a)$. The last condition in \mref{alpha} implies that
\[ \exists x\in R : a\md{x^{2}}4. \]
We have $(\frac{-d}a)=1$, since \A\
is not ramified at $a$, and hence
\[ \exists m\in R : -d\md{m^{2}}a. \]

Now, if 
\begin{equation}
  \label{eone}
  e_{1}=\frac{x+i}2 \mbox{ and } e_{2}=\frac{mi+ij}a,
\end{equation}
then the norms and traces of $e_{1},\,e_{2}$ and their products belong
to $R$, since
\[
\begin{array}{l}
  N(e_{1})=\frac{x^{2}-a}4,\;Tr(e_{1})=x,\\
  N(e_{2})=-\frac{d+m^{2}}a,\;Tr(e_{2})=0 \mbox{ and } Tr(e_{1}e_{2})=m.
\end{array}
 \]
Hence, we get that $\OO=\left<1,e_{1},e_{2},e_{1}e_{2}\right>$ is an order. The
matrix which takes $1,i,j,ij$ to $1,e_{1},e_{2},e_{1}e_{2}$ has
determinant equal to $\frac1{4a}$. Hence, $[\OO : \OO^{'}]=(4a)$ and we
have proved the following:

\prop{ \label{orderbasis}
Let $\A$ be a quaternion algebra over $K$ with $d(\A)=(d)$
principal. Let $a$ be chosen according to $\mref{alpha}$, so that
$\A\cong(a,-d)_{K}$. Then
$\OO=\left<1,e_{1},e_{2},e_{1}e_{2}\right>$ is a maximal order in \A, where
$e_{1}$ and $e_{2}$ are defined by $\mref{eone}$.
}

We remind that in general there may be other non-isomorphic maximal
orders in \A. However, if \A\ is for example a totally indefinite quaternion
algebra over a field $K$ with class number of $K$ equal to $1$, then
the maximal order is unique up to isomorphism \cite{Eichler4}.

\bibliography{sj}

\end{document}